\documentclass[sans]{amsart}
\usepackage{amsmath, amssymb}

\usepackage[square,sort&compress,comma,numbers]{natbib}
\usepackage{graphicx}
\usepackage{enumerate}
\usepackage{hyperref, tikz, times, xr-hyper}

\newtheorem{theorem}{Theorem}[section]
\newtheorem{lemma}[theorem]{Lemma}
\newtheorem{corollary}[theorem]{Corollary}
\newtheorem{proposition}[theorem]{Proposit{i}on}

\theoremstyle{definition}
\newtheorem{definition}[theorem]{Definit{i}on}

\theoremstyle{remark}

\numberwithin{equation}{section}

\newcommand{\abs}[1]{\left\lvert #1 \right\rvert}

\newcommand{\C}{\ensuremath{\mathbb{C}}}
\newcommand{\eps}{\ensuremath{\varepsilon}} 
\newcommand{\G}{\ensuremath{\mathbb{Z} [i]}}

\newcommand{\N}{\ensuremath{\mathbb{N}}}
\renewcommand{\O}{\ensuremath{\mathcal{O}}}
\newcommand{\Q}{\ensuremath{\mathbb{Q}}}
\newcommand{\R}{\ensuremath{\mathbb{R}}} 
\newcommand{\SLn}[2]{\ensuremath{\mathrm{SL}_{#1} \left( #2 \right)}}
\newcommand{\SLtwoC}{\ensuremath{\mathrm{SL}_2 \left( \mathbb{C} \right)}}
\newcommand{\SLtwoG}{\ensuremath{\mathrm{SL}_2 \left( \G \right)}}
\newcommand{\SLtwoOK}{\ensuremath{\mathrm{SL}_2 \left( \O_K \right)}}
\newcommand{\SLtwoZ}{\ensuremath{\mathrm{SL}_2 \left( \mathbb{Z} \right)}}

\newcommand{\Z}{\ensuremath{\mathbb{Z}}}

\begin{document}

  \title[Diophant{i}ne exponents for linear act{i}on of $\mathrm{SL}_2$ over discrete rings]{\protect{Diophant{i}ne exponents for standard linear act{i}ons of~$\mathbf{SL_2}$~over discrete rings in \C}}

  \author{L.~Singhal}
  \address{School of Mathemat{i}cs, Tata Inst{i}tute of Fundamental Research, Mumbai 400\,005}

  \email{singhal\symbol{64}math.t{i}fr.res.in}
\thanks{The author thanks his supervisor Dr.~Anish Ghosh for sugges{t}{i}ng the problem, reading the various preliminary versions of the manuscript and for being a constant source of encouragement. Financial support from CSIR, Govt.\ of India under SPM-07/858(0199)/2014-EMR-I is duly acknowledged.}

  \date{\today}


\begin{abstract}
  We give upper and lower bounds for various Diophant{i}ne exponents associated with the standard linear act{i}ons of \SLtwoOK\ on the punctured complex plane $\C^2 \setminus \{ \mathbf{0} \}$, where $K$ is a number f{i}eld whose ring of integers $\O_K$ is discrete and within a unit distance of any complex number. The results are similar to those of \citeauthor*{LN12} for \SLtwoZ\ act{i}on on $\R^2 \setminus \{ \mathbf{0} \}$ albeit for us, uniformly nice bounds are obtained only outside of a set of null measure.
\end{abstract}
\maketitle
  
  \section{Introduct{i}on}\label{S:intro}
    The set \Q\ of rat{i}onal numbers is dense in \R. However, one of the first works which tried to quant{i}fy this density came only in the nineteenth century from Dirichlet who stated that for any real number $\theta$ and all $Q > 1$, there exist integers $p$ and  $q, 1 \leq q < Q$ such that
    \begin{equation}
      \abs{\,q \theta - p\,} < \frac{1}{Q} < \frac{1}{q}.
    \end{equation}
    This is a consequence of the pigeon-hole principle (also known as Dirichlet's box principle). The inhomogeneous version was given by Minkowski using geometry of numbers. For any $\theta \in \R \setminus \Q$ and $\alpha \notin \Z \theta + \Z$, there exist integers $p$ and $q$ for which
    \begin{equation}
      \abs{\,q \theta - \alpha - p\,} < \frac{1}{4 \abs{q}}.
    \end{equation}
    Other than the fact that this second statement is only true for irrat{i}onal $\theta$, the error est{i}mate is also weak here than in Dirichlet's theorem where it is in terms of $Q^{-1} < q^{-1} (q < Q)$. If we now take two such inequalit{i}es with dif{f}erent $\alpha$'s, we are in the realm of simultaneous inhomogeneous Diophant{i}ne approximat{i}on~\citep{Cas}. Otherwise said, we are looking for inf{i}nitely many integral solut{i}ons $( p_1, q_1, p_2, q_2 )$ to the system of inequalit{i}es
    \begin{equation}
      \max \{ \abs{ q_1 \theta - p_1 - \alpha_1}, \abs{ q_2 \theta - p_2 - \alpha_2} \} < \eps.
    \end{equation}
    An extra demand that the pairs $( q_1, p_1 )$ and $( q_2, p_2 )$ be primit{i}ve can be fulf{i}lled by asking that the matrix
    \begin{equation}
      \begin{pmatrix}
        q_1 & p_1\\
        q_2 & p_2
      \end{pmatrix} \in \SLtwoZ.
    \end{equation}
    In recent t{i}mes, mathemat{i}cians have been interested in understanding more generally the nature of dense orbits for the act{i}on of a group $G$ on a homogeneous space $X$. In this respect, see the works of \citeauthor*{GGN3} \citep{GGN3, GGN4} where they relate the rate of approximation by `rat{i}onal points' on a homogeneous space $X$ of a semisimple group $G$ to the automorphic representat{i}ons of $G$ and compute the exact exponents for a number of examples. We point out upfront that their exponent $\kappa_{\Gamma}$ is exactly the inverse of the value $\hat{\mu}_{\Gamma}$ introduced in Def.~\ref{D:dexpunf} below.\par

    Let $K$ be any number f{i}eld whose ring of integers $\mathcal{O}_K$ is a discrete subring of \C. In addit{i}on, we require that any complex number $z$ should be within unit distance of some element of $\O_K$. The only such rings correspond to the rings of integers for the quadrat{i}c number f{i}elds $\Q (\sqrt{-d})$ where $d = 1, 3, 7$ or $11$ \cite[see][Remark~2.4]{Dan15}. By $\Gamma = \SLtwoOK$, we shall mean the lat{t}ice consist{i}ng of special linear matrices with entries belonging to $\O_K$. Consider its act{i}on on the punctured complex plane $\C^2 \setminus \{ 0 \}$ via matrix mult{i}plicat{i}on on the left. Abusing notat{i}on, we use $\abs{\cdot}$ with matrices as well as complex numbers which means that some clarif{i}cat{i}on is in order. For any matrix $A$, we let $\abs{A}$ be the maximum of the modulus of its entries while for a complex number, $\abs{z}$ stands for the Euclidean distance to the origin. We use lowercase Greek and both upper and lowercase Roman letters for various operat{i}ng matrices and vectors will be in boldface (e.~g.\ $\mathbf{z}$). The following terminology is mot{i}vated from \citeauthor{BL05}~\citep{BL05}.
    \begin{definition}\label{D:dexp}
      Let $\mathbf{z}, \mathbf{y} \in \C^2$. The Diophant{i}ne exponent $\mu_{\Gamma} ( \mathbf{z}, \mathbf{y} )$ is the quant{i}ty
      \[
         \sup \left\{ \omega \mid \abs{\gamma \mathbf{z} - \mathbf{y}} \leq | \gamma |^{- \omega} \textrm{ has inf{i}nitely many solut{i}ons in } \gamma \in \Gamma \right\}.
      \]
    \end{definition}
    \begin{definition}\label{D:dexpunf}
      The exponent $\hat{\mu}_{\Gamma} ( \mathbf{z}, \mathbf{y} )$ stands for the supremum of all $\omega$'s for which the system of inequalit{i}es
      \[
        \abs{\gamma \mathbf{z} - \mathbf{y}} \leq T^{- \omega},\quad | \gamma | \leq T
      \]
      has solut{i}ons for \emph{all} $T$ suf{f}iciently large.
    \end{definition}
    If the results of Dirichlet and Minkowski were to be recast in this language, they will say that the uniform exponents for respect{i}vely approximat{i}ng the points $( \theta, 0 )$ and $(\theta, \alpha)$ using an integral pair $(p, q)$ are both $\geq 1$. Further, measure\,-\,theoret{i}c considerat{i}ons dictate that both the equalit{i}es hold except on some set of Lebesgue measure zero. Prop.~\ref{P:minkowski} and the subsequent discussion gives us an analogue of Minkowski's theorem for approximat{i}ng a complex pair $(\xi, z)$ with the help of $\O_K$\,-\,integers.\par
    It follows from the def{i}nit{i}ons that for all pairs $\mathbf{z}, \mathbf{y} \in \C^2$, we have $\mu_{\Gamma} ( \mathbf{z}, \mathbf{y} ) \geq \hat{\mu}_{\Gamma} ( \mathbf{z}, \mathbf{y} )$. For the analogous situat{i}on of \SLtwoZ\ act{i}ng on $\R^2$, \citeauthor{LN12} \cite{LN12} came up with est{i}mates for the exponents def{i}ned above and in some cases, get the lower and upper bounds to be equal, which turn out to be funct{i}ons of the irrat{i}onality measures of the slopes of the starting and the target point. For approximat{i}on in our set{t}{i}ng, we follow in their footsteps to a large extent. When $d = 1$ so that $\O_K$ is the ring of Gaussian integers, a cont{i}nued fract{i}on expansion algorithm for complex numbers with partial quot{i}ents from $\O_K$ was given by  \citeauthor{Hur}~\cite{Hur}. The lat{t}er is made use of for construct{i}ng certain \emph{convergent matrices}. The case $d = 3$ has the ring of Eisenstein integers as its integral ring and we have an analogue in the shape of nearest (Eisenstein) integer algorithm~\cite{Dan15}. These help us to approach any fixed target point $\mathbf{y} \in \C^2$ start{i}ng at some \emph{``irrat{i}onal''} vector $\mathbf{z} = ( z_1, z_2 )^t \in \C^2$ both of whose coordinates are non-zero and the slope $\xi = z_1 / z_2 \in \C' := \C \setminus K$. Note that $\C'$ is a full measure subset of \C\ as $K$ is only countable.\par
    We give the following general def{i}nit{i}on inspired from that of irrat{i}onality measure for real irrat{i}onal numbers.
    \begin{definition}\label{D:irrmeasure}
      The $\mathbf K$\,\textbf{-\,irrat{i}onality measure} $\omega_K (z)$ for any $z \in \C'$ is the supremum of all numbers $\omega$ such that the inequality
      \[
        | q z - p | \leq \frac{1}{| q |^\omega}
      \]
      has inf{i}nitely many solut{i}ons in $p \in \O_K,\ q \in \O_K \setminus \{ 0 \}$.
    \end{definition}
    \citeauthor{Sul}~\cite[Theorem~1]{Sul} amongst others has formulated and proved the Khintchine theorem for Diophant{i}ne approximat{i}on of complex numbers by rat{i}onals coming from some fixed imaginary quadrat{i}c extension of \Q. In par{t}{i}cular, it implies that for all f{i}elds $K$ being considered here, the irrat{i}onality measure $\omega_K (z)$ is an almost everywhere constant funct{i}on on $\C'$ with respect to the induced Lebesgue measure. Using the convergence case of Borel-Cantelli lemma along with Dirichlet's box principle (see \citep[pg.~1]{Cas} and also Lem.~\ref{L:Dirichlet} below), one can independently verify its generic value to be $1$ and greater than that everywhere else. At this point, we remind the reader that the exponents $\mu_{\Gamma}$ and $\hat{\mu}_{\Gamma}$ def{i}ned above are invariant under $\Gamma$\,-\,act{i}on and, therefore, constant a.~e.\ owing to the ergodicity of the act{i}on.
    
    For non\,-\,negat{i}ve functions $f$ and $g$, the Vinogradov notat{i}on $f \ll g$ (similarly $f \gg g$) means that there exists some $C > 0$ for which $f (x) \leq C g (x)$ for all $x$ in the domain. The dependence of this implicit constant on some ambient parameters $a, b, c, \ldots$ will be o{f}{t}en indicated in the subscript as $\ll_{a, b, c, \ldots}$ The main result contained in this paper is given below.
    \begin{theorem}\label{Th:main}
      Let $K$ be of the form $\Q (\sqrt{-d})$ where $d = 1, 3, 7$ or $11$. Also, suppose that a cont{i}nued fract{i}on algorithm for approximat{i}ng an arbitrary complex number $z$ with elements of $K$ exists and has the following propert{i}es for all $n \gg 0$,
      \begin{enumerate}
        \item the denominators of the convergents rise monotonically, i.~e., $| q_{n + 1} | > | q_n |$, and
        \item there exists $r_0 \in \N$ and $\theta > 1$ for which $| q_{n + r_0} | \geq \theta | q_n |$.\label{I:theta}
      \end{enumerate}
      Then, for the full measure subset $\{ \mathbf{z} = ( z_1, z_2 )^t \in \C^2 \mid z_1 / z_2 \in \C \setminus K,\ \omega_K ( z_1 / z_2 ) = 1 \}$ and $\Gamma = \SLtwoOK$ act{i}ng linearly on the complex plane, it is true that:
      \begin{enumerate}[i)]
      \item the exponent of approximation to the origin is $\mu_{\Gamma} ( \mathbf{z}, \mathbf{0} ) = \hat{\mu}_{\Gamma} ( \mathbf{z}, \mathbf{0} ) = 1$,\label{I:origin}
      \item for almost all target points $\mathbf{y}$ with the slope $y = y_1 / y_2 \in \C \setminus K$ and $\omega_K (y) = 1$,\label{I:irr}
      \begin{equation}\label{E:mirr}
        1/3 \leq \hat{\mu}_{\Gamma} ( \mathbf{z}, \mathbf{y} ) \leq \mu_{\Gamma} ( \mathbf{z}, \mathbf{y} ) \leq 1/2, \textrm{ and }
      \end{equation}
      \item for target points $\mathbf{y}$ with slope $y \in K$, we have\label{I:rat}
      \begin{equation}\label{E:mrat}
        \mu_{\Gamma} ( \mathbf{z}, \mathbf{y} ) = \hat{\mu}_{\Gamma} ( \mathbf{z}, \mathbf{y} ) = 1/2.
      \end{equation}
      \end{enumerate}
    \end{theorem}
    \noindent While the discreteness of $\O_K$ (which is ensured by taking $K = \Q (\sqrt{-d} ),\ d$ as above)  immediately implies that of \SLtwoOK\ and is also used in working out the generic $K$\,-\,irrat{i}onality measure $\omega_K$, the exponent{i}al rise of denominator sizes in the cont{i}nued fract{i}on algorithm helps in bounding the various intermediate matrices properly.\par
    
    We emphasize again that the results of \citeauthor{LN12}~\citep{LN12} are valid for all start{i}ng points $\mathbf{x} \in \R^2 \setminus \{ \mathbf{0} \}$ with irrat{i}onal slopes and all target points, but we have nice answers only on a full measure subset of star{t}{i}ng points where $\omega_K ( z_1 / z_2 ) = 1$. Furthermore, we get a reasonable lower bound for target points $\mathbf{y}$ with a $K$\,-\,irrat{i}onal slope only when $\omega_K ( y_1 / y_2 ) = 1$ and the upper bound of $1/2$ in Eq.~\eqref{E:mirr} is true for some full measure subset of $\C^2$ (coming out of the Borel\,-\,Cantelli lemma and perhaps depending on $\mathbf{z}$).\par
    
    A con{t}{i}nued frac{t}{i}on theory as assumed in Th.~\ref{Th:main} is provided for $d =1$ and $3$ in \cite{Hur, Dan15, DN14}. We discuss and suitably modify some of their statements in the next sect{i}on. The jugglery with approximat{i}ng matrices comes in Sect.~\ref{S:conmat} which, in various parts, gives us Theorem~\ref{Th:main} (Prop.~\ref{P:origin} for \eqref{I:origin}, Props.~\ref{P:irr} \& \ref{P:generic} for~\eqref{I:irr} and Prop.~\ref{P:rat} for~\eqref{I:rat}, respect{i}vely).

  \section{\protect{Cont{i}nued fract{i}ons for complex numbers}}\label{S:cfc}
    \citeauthor{DN14} \cite{DN14} have considered a family of cont{i}nued fract{i}on expansions for complex numbers where the partial quot{i}ents $a_n \in \Z [i]$, the ring of Gaussian integers. In~\cite{Dan15}, \citeauthor*{Dan15} also dealt with cont{i}nued fract{i}ons in terms of Eisenstein integers $a + b \zeta$ where $a, b \in \Z$ and $\zeta^2 + \zeta = -1$.
    
    In part{i}cular, \citeauthor{DN14} give the best known results for the rate of approximat{i}on by convergents coming from Hurwitz's algorithm. Hurwitz \cite{Hur} described a simple \emph{nearest integer algorithm} which picks a Gaussian integer $a$ nearest to any given complex number $z$ (if there is more than one candidate sat{i}sfying the condit{i}on, choose any one of them). One then proceeds by induct{i}on as
    \begin{equation}
      z_0 = z \textrm{ and } z_{n + 1} = ( z_n - a_n )^{-1}.
    \end{equation}
    If $z = [\,a_0, a_1, \ldots\,]$, then $a_i \in \G$ for all $i$ and $| a_i | > 1$ for $i \geq 1$. On def{i}ning the associated numerator and denominator (of the $n$th convergent) sequences of Gaussian integers in a recursive fashion as
    \begin{align}\label{E:Qpair}
       p_{-2} = 0,\quad &p_{-1} = 1,\quad&p_n &= a_n p_{n - 1} + p_{n - 2} \textrm{ for } n \geq 0, \textrm{ and }\notag\\
       q_{-2} = 1,\quad &q_{-1} = 0,\quad&q_n &= a_n q_{n - 1} + q_{n - 2} \textrm{ for } n \geq 0,
    \end{align}
    we are assured of the exponent{i}al growth of the size of the denominators, namely that $| q_n | > | q_{n - 1} |$ and $| q_n | \geq \theta | q_{n - 2} |$ for all $n \geq 1$, where $\theta = ( \sqrt{5} + 1 ) / 2$ \cite[Corollary~5.3]{DN14}. The lat{t}er guarantees that the distance between the complex number $z$ and its $n$th convergent is small enough in terms of the size of the denominator $q_{n + 1}$ of the succeeding convergent.
    
    We have a similar situat{i}on at hand for the Eisenstein integers. Theorem~4.3 of~\cite{Dan15} tells us that for the cont{i}nued fract{i}on expansion with respect to the nearest integer algorithm, we have the monotonous rise of the denominator sizes as well as $| q_n | > 4 | q_{n - 2} | / 3$ for $n \geq 1$. We now give an alterat{i}on of \cite[Proposit{i}on~2.1]{DN14}.
    \begin{lemma}\label{L:cfeub}
      Let $R$ be a discrete subring of \C\ with $\mathrm{Frac} (R)$ being its field of fract{i}ons. Further, let $\{ a_n \} \subset R$ be a sequence which def{i}nes a cont{i}nued fract{i}on expansion of some $z \in \C \setminus \mathrm{Frac} (R)$ and $p_n / q_n$ be the corresponding sequence of convergents for which the hypothesis of Th.~\ref{Th:main} holds. Then, there exist $C_1$ and $n_0$ posit{i}ve such that
      \begin{equation}\label{E:uerror}
        \left| q_n z - p_n \right| \leq \frac{C_1}{| q_{n + 1} |}\quad \forall n > n_0.
      \end{equation}
    \end{lemma}
    \begin{proof}
      We need to look more closely at the proof given in~\cite{DN14} which goes through for any discrete ring without any changes whatsoever. There, the authors have argued that
      \begin{equation}\label{E:refine}
         \left| z - \frac{p_n}{q_n} \right| \leq \sum_{k = 0}^\infty \frac{1}{| q_{n + k} q_{n + k + 1} |} \leq \frac{C_0}{| q_n |^2}, \textrm{ where } C_0 = \frac{r_0 \theta^2}{\theta^2 - 1}.
      \end{equation}
      We separate the first term from the series on the right above
      \begin{equation}
        \frac{1}{|q_n q_{n + 1} |} + \sum_{k = 1}^\infty \frac{1}{| q_{n + k} q_{n + k + 1} |} \leq \frac{1}{|q_n q_{n + 1} |} + \frac{C_0}{|q_{n + 1} |^2}
      \end{equation}
      and the upper bound is arrived at by taking $n = n + 1$ in the last step of their calculat{i}on. Mult{i}plying Eq.~\eqref{E:refine} by $| q_n |\ ( \neq 0 )$ on both sides and recalling that $| q_n | < | q_{n + 1} |$, we then have that the scaled error $| \epsilon_n | = | q_n z - p_n | \leq C_1 | q_{n + 1} |^{-1}$, where $C_1 = C_0 + 1$ is an absolute constant.
    \end{proof}
    
    \noindent Let us now try to obtain a lower bound for $\epsilon_n$ which is not a priori available. Unlike the simple cont{i}nued fract{i}ons for real numbers, we get a very weak lower est{i}mate for the $n$-th error term. But before that, consider two di{f}{f}erent convergents $p_n / q_n$ and $p_{n + r} / q_{n + r}$ for some $n \geq 0,\ r > 0$ arising from a cont{i}nued fract{i}on expansion of a f{i}xed $z \in \C \setminus K$, where the associated part{i}al quot{i}ents belong to the (discrete) ring of integers $\O_K$. Our claim is that the two convergents are not the same complex number. If not, let $p_{n + r} = \kappa p_n$ and $q_{n + r} = \kappa q_n$ for some $\kappa \in \C \setminus \{ 0 \}$. As $| q_{n + 1} | > | q_n |$ for all $n$, we get that $| \kappa | > 1$. Also,
    \begin{equation}
      \left| \kappa ( p_n q_{n + r - 1} - q_n p_{n + r - 1} ) \right| = \left| p_{n + r} q_{n + r - 1} - q_{n + r} p_{n + r - 1} \right| = 1
    \end{equation}
    and hence, the non-zero complex number $p_n q_{n + r - 1} - q_n p_{n + r - 1} \in \O_K$ has absolute value at least $1$, the lat{t}er being a discrete ring. But, this is a contradict{i}on.
    \begin{lemma}\label{L:cfelb}
      With the same notat{i}ons and convent{i}ons as in Lemma~\ref{L:cfeub}, there exist $C_2 > 0$ and $r_1 \in \N$ such that
      \[
        | \epsilon_n | \geq \frac{C_2}{| q_{n + r_1} |}.
      \]
    \end{lemma}
    \begin{proof}
      Apply triangle inequality to the three numbers $z, p_n / q_n$ and $p_{n + r} / q_{n + r}$ giving
      \begin{equation}
        \left| z - \frac{p_n}{q_n} \right| \geq \left| \frac{p_{n + r}}{q_{n + r}} - \frac{p_n}{q_n} \right| - \left| z - \frac{p_{n + r}}{q_{n + r}} \right| \geq \frac{1}{| q_n | | q_{n + r} |} - \frac{C_0}{| q_{n + r} |^2}
      \end{equation}
      when we employ~\eqref{E:refine}. This implies that
      \begin{equation}
        | \epsilon_n | = | q_n z - p_n | \geq \frac{1}{| q_{n + r} |} - \frac{| q_n |}{| q_{n + r} |}\frac{C_0}{| q_{n + r} |}.
      \end{equation}
      Now, $| q_{n + r} | \geq \theta\,| q_{n + r - r_0} | \geq \ldots \geq \theta^{\lfloor \frac{r}{r_0} \rfloor} | q_n |$ for all $n, r$ by our assumpt{i}on. Thus,
      \begin{equation}
        | \epsilon_n | \geq \frac{1}{| q_{n + r} |} - \frac{C_0}{\theta^{\lfloor \frac{r}{r_0} \rfloor} \left| q_{n + r} \right|}
      \end{equation}
      As $\theta > 1$, the constant in the second term on the right side becomes less than $1$ for some $r_1$ suf{f}{i}ciently large and we get the required lower bound for some constant $C_2 > 0$ and $r_1 \in \N$ depending only on $R$ and the cont{i}nued fract{i}on algorithm in ef{f}ect.
    \end{proof}
    When the $K$\,-\,irrat{i}onality measure $\omega_K (z)$ is f{i}nite and $\omega > \omega_K (z)$, then we must have $| q_{n + 1} | \leq | q_n |^{\omega}$ for all $n \geq N_0 ( \omega )$. Combining Lemmata~\ref{L:cfeub} and \ref{L:cfelb}, we conclude that for all $\omega > \omega_K (z)$,
    \begin{equation}\label{E:finalerror}
      \frac{C_2}{\left| q_{n + 1} \right|^{\omega^{r_1 - 1}}} \leq | \epsilon_n | \leq \frac{C_1}{\left| q_{n + 1} \right|}
    \end{equation}
    for all large enough $n$. In addit{i}on, we get the usual ident{i}ty
    \begin{equation}\label{E:coprimecf}
         q_n p_{n - 1} - p_n q_{n - 1} = ( -1 )^n
      \end{equation}
      as a bonus from the formal theory of cont{i}nued fract{i}ons. It is this part{i}cular property of theirs which enables them to be of good use in construct{i}ng the so-called \emph{convergent matrices} discussed in the next sect{i}on.

  \section{Convergent matrices}\label{S:conmat}
    Let $\xi \in \C \setminus K$ and let $p_k / q_k$ for $p_k, q_k \in \O_K$ denote the convergent of order $k$ to $\xi$, due to some cont{i}nued fract{i}on expansion algorithm which sat{i}sf{i}es the hypothesis of Theorem~\ref{Th:main}. The construct{i}on of convergent matrices for the complex set{t}{i}ng mimics the one for $\R^2$. As in~\cite{LN12}, we def{i}ne the $k$-th convergent matrix
    \begin{equation}\label{E:conmat}
      M_k := \begin{pmatrix}
               q_k                   & -p_k\\
	       (-1)^{k - 1}q_{k - 1} & (-1)^k p_{k - 1}
	     \end{pmatrix} \in \SLtwoOK\ = \Gamma.
    \end{equation}
    The powers of $-1$ have been inserted so that the matrices are special linear ones once we have~\eqref{E:coprimecf}. The supremum norm of the above matrix is $\max ( | q_k |, | p_k | )$ since the size of the denominators increases monotonically, and the numerators $p_k$'s should increase accordingly in order to approximate better and better the fixed complex number $\xi$. If necessary, we pre-mult{i}ply the vector $\mathbf{z} \in \C^2$ by the \SLtwoOK\ matrix
    \begin{equation}\label{E:J}
      J = \begin{pmatrix}
             0 & -1\\
             1 & 0
          \end{pmatrix}
    \end{equation}
   to have the slope $\xi$ with $| \xi | \leq 1$ while the size $| J \mathbf{z} | = | \mathbf{z} |$ remains the same. Also, note that $| J \gamma | = | \gamma J | = \abs{\gamma}$ for any $2 \times 2$ matrix $\gamma$. Thus, it is alright to take $| M_k | \asymp | q_k |$ and we will do as much from here on without explicitly saying so. When concerned with the $\Gamma$\,-\,orbit of the point $\mathbf{z} = ( z_1, z_2 )^t$ having $z_1 / z_2 = \xi$, we see that
    \begin{equation}
      M_k \mathbf{z} = z_2
                       \begin{pmatrix}
                         \epsilon_k\\
	                 (-1)^{k - 1} \epsilon_{k - 1}
	     \end{pmatrix}
    \end{equation}
    implying that
    \begin{equation}\label{E:lbmu}
      | M_k \mathbf{z} | \leq | \mathbf{z} |\frac{C_1}{| q_k |} \ll_K \frac{ | \mathbf{z} |}{| M_k |} \textrm{ as both } | \epsilon_k | \leq \frac{C_1}{| q_{k + 1} |} < \frac{C_1}{| q_k |} \textrm{ and } | \epsilon_{k - 1} | \leq \frac{C_1}{| q_k |}
    \end{equation}
    leveraging~\eqref{E:uerror}. Thus, there are infinitely many such matrices and this immediately tells us that $\mu ( \mathbf{z}, \mathbf{0} ) \geq 1$. The proof of $\mu ( \mathbf{z}, \mathbf{0} ) \leq 1$ goes along the same lines as~\cite[Lemma~1]{LN12} with the necessary modif{i}cat{i}ons by  $| q_k |$ and $| q_{k + 1} |$ replacing $q_k$ and $q_{k + 1}$, respect{i}vely and having appropriate constants in place. Trivially, we also get an upper bound on $\hat{\mu} ( \mathbf{z}, \mathbf{0} )$. For the reverse inequality, it suf{f}{i}ces to consider the matrices $M_k$ with
    \begin{align}
      \left| q_k \right| \asymp | M_k | \leq T &\leq | M_{k + 1} | \asymp \left| q_{k + 1} \right|, \textrm{ and}\notag\\
      \left| M_k \mathbf{z} \right| \ll_{\mathbf{z}} \frac{1}{\left| q_k \right|} &\ll \frac{1}{\left| q_{k + 1} \right|^{1 / \omega}} \ll \frac{1}{T^{1 / \omega}}
    \end{align}
    for $\omega > \omega_K (\xi)$ and all $k > k_0 = k_0 ( \xi, \omega)$. Let{t}{i}ng $\omega_K (\xi) \leftarrow \omega$ from the right, we have
   \begin{proposition}\label{P:origin}
     For any vector $\mathbf{z} = ( z_1, z_2 )^t \in \C^2$ with slope $\xi = z_1 / z_2 \in \C'$ such that a cont{i}nued fract{i}on expansion for $\xi$ in terms of $\O_K$\,-\,integers exist and satisf{i}es the condit{i}ons of Th.~\ref{Th:main}, we have the exponents of approximat{i}on
     \[
       \mu_{\Gamma} ( \mathbf{z}, \mathbf{0} ) = 1\quad \textrm{and}\quad \frac{1}{\omega_K (\xi)} \leq \hat{\mu}_{\Gamma} ( \mathbf{z}, \mathbf{0} ) \leq 1.
     \]
   \end{proposition}
   \noindent We digress now for a bit to prove a claim we made in the discussion after Def.~\ref{D:irrmeasure}.
   \begin{lemma}\label{L:Dirichlet}
     For $K = \Q ( \sqrt{-d} )$ where $d = 1, 3, 7$ or $11$, the $K$\,-\,irrat{i}onality measure $\omega_K (z)$ is equal to $1$ for Lebesgue almost all and $\geq 1$ for all $z \in \C'$.
   \end{lemma}
   \begin{proof}
     Given $Q > 1$, the number of $\O_K$ integers $q$ with $\abs{q} \leq Q/2$ is $\geq c_K Q^2$ for some $c_K > 0$. For each of these $q$'s, there exists a unique $p = p(q) \in \O_K$ such that the complex number $q z - p$ belongs to a fixed fundamental polygon $\mathcal{F}_K$ for $\O_K$ in \C. Therefore, we have at least $c_K Q^2$ many distinct numbers in $\mathcal{F}_K$ as $z \notin K$. If we now divide this polygon into $\approx c_K Q^2$ many subpolygons each of which has diameter $\leq c^{-1/2}_K Q^{-1}$, by Dirichlet's pigeonhole principle, we should have that one of them contains both $q_1 z - p_1$ and $q_2 z - p_2$ for some $\abs{q_1}, \abs{q_2} \leq Q/2$ and $q_1 \neq q_2$. In conclusion,
     \begin{equation}\label{E:KDirichlet}
       \abs{ ( q_1 - q_2 ) z - ( p_1 - p_2 ) } \leq \frac{1}{\sqrt{c_K} Q} \ll_K \frac{1}{| q_1 - q_2 |}
     \end{equation}
     giving us that $\omega_K (z) \geq 1$. To see that the equality holds for almost all $z \in \C'$, notice the number of $p \in \O_K$ such that $p / q \in \mathcal{F}_K$ is $\leq b_K \abs{q}^2$ for any fixed $q$ and the number of $q \in \O_K$ for which $\abs{q} \sim Q$ is $\leq b'_K Q$. Therefore, the series in the Borel\,-\,Cantelli lemma for the family of discs of radius $1 / \abs{q}^{1 + s}$ around the $K$\,-\,rat{i}onal point $p/q$ is dominated by
     \begin{equation}\label{E:bcseries}
       \sum_{Q > 1} b_K b'_K \frac{Q^3}{(Q^{1 + s})^2}.
     \end{equation}
     The la{t}{t}er converges for all $s > 1$ implying that for $s$ in this range, the limsup set
     \begin{equation}\label{E:limsupset}
       \limsup_{p/q \in K} D \left( \frac{p}{q}, \frac{1}{\abs{q}^{1 + s}} \right)
     \end{equation} has Lebesgue measure zero. In other words, $\omega_K (z) = 1$ for almost all $z \in \C'$.
   \end{proof}
   \noindent Hence, the generic value (in $\mathbf{z}$) of both $\mu_{\Gamma} ( \mathbf{z}, \mathbf{0} )$ and $\hat{\mu}_{\Gamma} ( \mathbf{z}, \mathbf{0} )$ is $1$. We, thereby, have the f{i}rst claim of Theorem~\ref{Th:main}.\par
   
   A func{t}{i}on $h : X \rightarrow [ 0, \infty )$ on a countable space $X$ is said to be a \emph{height func{t}{i}on} if for each $Q \geq 0$, the set $h^{-1} [ 0, Q ]$ is f{i}nite.
   If $( \varphi, G )$ is an ac{t}{i}on of a countable group $G$ with height func{t}{i}on $h$ on a metric space $X$, the exponent $\mu_{\varphi} ( x, y)$ stands for
   \begin{equation}\label{E:muphi}
     \sup \{ \mu \mid \mathrm{dist} ( gx, y ) < h (g)^{- \mu} \textrm{ has inf.\ many solut{i}ons in } g \}.
   \end{equation}
   The uniform variant $\hat{\mu}_{\varphi} ( x, y )$ is given in the same fashion for all $x, y \in X$. Next, we make a simple and more general observat{i}on whose proof is immediate from the de{f}{i}ni{t}{i}ons.
    \begin{proposition}\label{P:Ghom}
      Let $( G_1, h_1 )$ and $( G_2, h_2 )$ be countable groups with $h_i$ being a height funct{i}on on $G_i$ and $\rho : G_1 \rightarrow G_2$ be a group homomorphism which respects $h_1$ in the sense that there exists $c > 1$ s.~t.
      \begin{equation}
        \frac{1}{c} h_1 (g) \leq h_2 ( \rho (g) ) \leq c h_1 (g)\ \forall g \in G_1.
      \end{equation}
      Further, let $( \varphi_i, G_i ),\ i = 1, 2$ be group act{i}ons on a metric space $X$ and $\varphi_2 \circ \rho = \varphi_1$. Then, for all pairs $x, y \in X$,
      \begin{equation}
        \mu_{\varphi_2} ( x, y ) \geq \mu_{\varphi_1} ( x, y ) \textrm{ and } \hat{\mu}_{\varphi_2} ( x, y ) \geq \hat{\mu}_{\varphi_1} ( x, y ).
      \end{equation}
    \end{proposition}
    \noindent As the group \SLtwoG\ sits inside \SLn{4}{\Z}\ owing to $a + i b \leftrightarrow \begin{pmatrix} a &-b\\ b &a \end{pmatrix}$ with the height funct{i}on on \SLtwoG\ preserved and the standard linear act{i}ons on $\C^2 \cong \R^4$ coinciding under the resultant embedding, we have the following corollary for simultaneous approximat{i}on by primit{i}ve integral vectors in dimension $4$ by combining Props.~\ref{P:origin} and \ref{P:Ghom}.
   \begin{corollary}\label{C:R4}
     For almost all $\mathbf{v} \in \R^4$, the exponents $\mu ( \mathbf{v}, \mathbf{0}_4 )$ and $\hat{\mu} ( \mathbf{v}, \mathbf{0}_4 )$ for approaching the origin via \SLn{4}{\Z} orbit, are greater than or equal to $1$.
   \end{corollary}
   \noindent We will not write down the corresponding true statements for other target points in $\R^4$. The next lemma bounds the size of a convergent matrix $\gamma \in \SLtwoOK$ in terms of the entries in its decomposit{i}on. The idea here is to bring the start{i}ng point $\mathbf{z}$ suf{f}iciently close to the origin using matrices $M_k$ as above, spread it around as a latt{i}ce with the help of the subgroup
   \begin{equation}
     \mathcal{U} = \left\{ U^{\ell} :=
                   \begin{pmatrix}
                      1 & \ell\\
                        & 1
                   \end{pmatrix} \big|\ \ell \in \O_K \right\},
   \end{equation}
   and finally rotate the lat{t}ice so obtained by applying matrices $N_j$ which attempt to take the ``complex line'' $\langle\,( z_1 , 0 )\,\rangle$ closer to $\langle\,( z_1 , \xi z_1 )\,\rangle$.
   \begin{lemma}[\citeauthor{LN12} \cite{LN12}]\label{L:gammabound}
      Let $k \in \N$ and $\ell \in \O_K$. For any arbitrary
      \begin{equation}
        N = \begin{pmatrix} t & t'\\ s & s' \end{pmatrix} \in \Gamma,
      \end{equation}
      the matrix $\gamma = N U^{\ell} M_k$ sat{i}sf{i}es
      \begin{equation}
         \left| \left( \ell q_{k - 1} + (-1)^{k - 1} q_k \right) s \right| - \left| s' q_{k - 1} \right| \leq \abs{\gamma} \ll \left| \ell q_{k - 1} \right| \left| N \right| + \left| N \right| \left| q_k \right|.
      \end{equation}
   \end{lemma}
   \begin{proof}
     After two matrix mult{i}plicat{i}ons
     \begin{align}
       \gamma &= N U^{\ell} M_k = \begin{pmatrix} t   & t'   \\ s                     & s'               \end{pmatrix}
                                  \begin{pmatrix} 1   & \ell \\ 0                     & 1                \end{pmatrix}
                                  \begin{pmatrix} q_k & -p_k \\ (-)^{k - 1} q_{k - 1} & (-1)^k p_{k - 1} \end{pmatrix}\notag\\
              &= \begin{pmatrix}
                   tq_k + (-1)^{k - 1} q_{k - 1} ( t  \ell + t') & - tp_k + (-1)^k p_{k - 1} (t \ell + t')\\
                   sq_k + (-1)^{k - 1} q_{k - 1} ( s \ell + s') & - sp_k + (-1)^k p_{k - 1} (s \ell + s')
                 \end{pmatrix}.
     \end{align}
     The bottom left entry of the matrix determines the lower bound in the lemma as soon as we employ the triangle inequality. Because we have already reduced to the case $| \xi | \leq 1$, for all large enough $n$ we have $| p_n | \ll | q_n |$ and then, the upper bound is easy enough.
   \end{proof}
   We now take steps towards obtaining bounds for the vector $\gamma ( \xi, 1 )^t$. The lemma below is again due to \citeauthor{LN12}. We sketch its proof here to merely point out the minor dif{f}erence(s) with the real case.
   \begin{lemma}\label{L:gammay}
     Let $k, \ell, N$ and $\gamma = N U^{\ell} M_k = \begin{pmatrix} v_1 & u_1\\ v_2 & u_2 \end{pmatrix}$ be as in Lemma~\ref{L:gammabound} and $y \in \C$. If $\delta = sy - t$ and $\delta' = s' y - t'$, we have that
     \[
       \left| v_1 \xi + u_1 - y \left( v_2 \xi + u_2 \right) \right| \ll_K \frac{| \delta \ell + \delta' |}{| q_k |} + \frac{| \delta |}{| q_{k + 1} |}.
     \]
   \end{lemma}
   \begin{proof}
     To begin with,
     \begin{align}
       y \left( v_2 \xi + u_2 \right) - \left( v_1 \xi + u_1 \right) &= ( -1\ \ y ) \gamma
       \begin{pmatrix}
         \xi\\
         1
       \end{pmatrix} = ( -1\ \ y )
       \begin{pmatrix}
         t & t'\\
         s & s'
       \end{pmatrix} U^{\ell} M_k 
       \begin{pmatrix}
         \xi\\
         1
       \end{pmatrix}\notag\\
                       &= ( \delta\quad \delta' ) U^{\ell} 
       \begin{pmatrix}
         \epsilon_k\\
         (-1)^{k - 1} \epsilon_{k - 1}
       \end{pmatrix}\\
                       &= \delta \epsilon_k + (-1)^{k - 1} ( \delta \ell + \delta' ) \epsilon_{k - 1},\notag
     \end{align}
     Since for the cont{i}nued fract{i}on expansions being studied here, we have $| \epsilon_n | \ll_K | q_{n + 1} |^{-1}$ for all $n \gg 1$, the claim follows.
   \end{proof}
   \noindent If $( \Lambda_1, \Lambda_2 )^t$ is the dif{f}erence $\gamma \mathbf{z} - \mathbf{y}$, then
   \begin{equation}\label{E:lambdas}
     \Lambda_1 = z_2 ( v_1 \xi + u_1 ) - y_1,\quad \Lambda_2 = z_2 ( v_2 \xi + u_2 ) - y_2.
   \end{equation}
   and on further choosing $y = y_1 / y_2$, we get
   \begin{align}\label{E:lambdadiff}
     \abs{\Lambda_1 - y \Lambda_2} &= \abs{z_2 \left( ( v_1 \xi + u_1 ) - y ( v_2 \xi + u_2 ) \right)}\notag\\
                                   &\ll_K \abs{z_2} \left( \frac{| \delta \ell + \delta' |}{| q_k |} + \frac{| \delta |}{| q_{k + 1} |} \right).
   \end{align}
   Once we bound one of the components (say $\Lambda_2$) and the difference $| \Lambda_1 - y \Lambda_2 |$, the vector $( \Lambda_1, \Lambda_2 )^t$ is bounded automat{i}cally. We proceed to do just that. Af{t}er a slight adjustment in the proof of Lemma~\ref{L:gammay}, we deduce
   \begin{align}\label{E:lambda2}
     \Lambda_2 &= z_2 ( v_2 \xi + u_2 ) - y_2\\
               &= z_2 \left( s \epsilon_k + (-1)^{k - 1} ( s \ell + s' ) \epsilon_{k - 1} \right) - y_2 = (-1)^{k - 1} z_2 s \epsilon_{k - 1} ( \ell - \rho ),\notag
   \end{align}
   where
   \begin{equation}\label{E:rho}
     \rho = \frac{(-1)^{k - 1} y_2}{z_2 s \epsilon_{k - 1}} - \frac{(-1)^{k - 1} \epsilon_k}{\epsilon_{k - 1}} - \frac{s'}{s}
   \end{equation}
   helps us to decide the value of $\ell$ such that $| \ell - \rho | \leq C_3$ for some constant $C_3$ depending only on $\O_K$ (or $K$ if you will) and $| \ell | \leq | \rho |$, having fixed $M_k$ and $N$ first.

  \subsection{Generic target points}\label{SS:generic}
    In the next few pages, we discuss the situat{i}on where the target point $\mathbf{y} = ( y_1, y_2 )^t \in \C^2$ has slope $y = y_1 / y_2 \in \C' = \C \setminus K $. As such points const{i}tute a set of full measure in $\C^2$, we shall be inferring propert{i}es of almost all points in the complex plane. Let $t_j / s_j$ and $t_{j - 1} / s_{j - 1}$ be consecut{i}ve convergents in an $\O_K$-cont{i}nued fract{i}on expansion of $y$ for our fixed target point $\mathbf{y}$. As argued for $\xi$, we may also suppose $|y| \leq 1$ thanks to the $J$ of Eq.~\eqref{E:J}. When
    \begin{equation}
      t = t_j,\ s = s_j,\ t' = (-1)^{j - 1} t_{j - 1} \textrm{ and } s' = (-1)^{j - 1} s_{j - 1},
    \end{equation}
    the matrix $N_j$ given by $\begin{pmatrix} t & t'\\ s & s' \end{pmatrix}$ belongs to $\Gamma = \SLtwoOK$ and for any $\omega > \omega_K (\xi)$, the auxiliary term $\rho$ in Eq.~\eqref{E:rho} is conf{i}ned within the range
    \begin{equation}
      \frac{1}{C_1} \abs{\frac{y_2 q_k}{z_2 s_j}} - 2 \leq \abs{\rho} \leq \frac{\abs{y_2} \abs{q_k}^{\omega^{r_1 - 1}}}{\abs{C_2 z_2 s_j}} + 2
    \end{equation}
    for all $k$ large enough using Lemmata~\ref{L:cfeub} and \ref{L:cfelb}, the fact that the error term
    \begin{equation}
      \epsilon_n = \frac{(-1)^n}{z_1 \cdots z_{n + 1}}
    \end{equation}
    with $| z_i | \geq 1$ for $i > 0$~\cite[Prop.~2.1 (i)]{Dan15} and a monotonous increase in the size of denominators $s_j$'s. This in turn tells us that the opt{i}mal choice of $\ell$ obeys
    \begin{equation}\label{E:ellrange}
      \frac{1}{C_1}\abs{\frac{y_2 q_k}{z_2 s_j}} - ( C_3 + 2 ) \leq \abs{\ell}\ \ll_{\mathbf{y}, \mathbf{z}, K} \frac{\abs{q_k}^{\omega^{r_1 - 1}}}{\abs{s_j}}  + 1.
    \end{equation}
    Subst{i}tut{i}ng this in Lemmata~\ref{L:gammabound} and \ref{L:gammay}, we have an $\O_K$\,-\,analogue of \cite[Lemma~4]{LN12}.
    \begin{lemma}\label{L:gammairr}
      Let $j \in \N, k \gg 0$ and $\omega > \omega_K (\xi)$. There exists $\gamma = N_j U^{\ell} M_k \in \Gamma$ for some $\ell \in \O_K$ such that
      \[
        \left| \frac{\abs{y_2}}{\abs{C_1 z_2}}\abs{q_k q_{k - 1}} - | s_j q_k | \right| -  ( C_3+3 ) | s_j q_{k - 1} | \leq | \gamma | \ll_{\mathbf{y}, \mathbf{z}, K} \abs{q_{k - 1}} | q_k |^{\omega^{r_1 - 1}} + | s_j q_k |
      \]
      as well as
      \begin{equation}
        | \gamma \mathbf{z} - \mathbf{y} | \ll_{\mathbf{y}, \mathbf{z}, K} \frac{| q_k |^{\omega^{r_1 - 1} - 1}}{\abs{s_j s_{j + 1}}} + \left| \frac{s_j}{q_k} \right|.
      \end{equation}
    \end{lemma}
    \begin{proof}
      The bounds for $\gamma$ are straightforward. Insofar as $\gamma \mathbf{z} - \mathbf{y}$ is concerned,
      \begin{equation}\label{E:lambda2bound}
        | \Lambda_2 | \leq C_1 C_3 \left| \frac{z_2 s_j}{q_k} \right|,
      \end{equation}
      using Eq.~\eqref{E:lambda2} and an opt{i}mal choice of $\ell$ as explained immediately after Eq.~\eqref{E:rho}. Moreover, $| \Lambda_1 | \leq | y \Lambda_2 | + | \Lambda_1 - y \Lambda_2 | \leq | \Lambda_2 | + | \Lambda_1 - y \Lambda_2 |$ as $| y | \leq 1$. The quant{i}ty $| \Lambda_1 - y \Lambda_2 |$ is bounded using Eq.~\eqref{E:lambdadiff} while we recall that $| \delta | \leq C_1 / | s_{j + 1} |$ and $| \delta' | \leq C_1 / | s_j |$.
    \end{proof}
    With the above lemma on our side, we now make an appropriate choice of the indices $j$ and $k$ so that
    \begin{equation}\label{E:jandk}
      | q_{k - 1} |^{1/3} < | s_j | \leq | q_k |^{1/3} < | s_{j + 1} |.
    \end{equation}
    We are assured of the existence of arbitrarily large pairs $(j, k)$ sat{i}sfying the inequalit{i}es~\eqref{E:jandk} as $| q_k |$'s and $| s_j |$'s are strictly increasing sequences of real numbers. Such a pair is then fed into the statement of Lemma~\ref{L:gammairr} to give us
    \begin{equation}\label{E:errorirrfinal}
      | \gamma \mathbf{z} - \mathbf{y} | \ll_{\mathbf{y}, \mathbf{z}, K} \frac{1}{| q_{k - 1} |^{\frac{1}{3}} | q_k |^{\frac{4}{3} - \omega^{r_1 -1}}} + \frac{1}{| q_k |^{\frac{2}{3}}}.
    \end{equation}
    In this work, we are mostly concerned with the $\Gamma$-orbits of generic points in $\C^2$ whose slope has $K$\,-\,irrat{i}onality measure equal to $1$. Thus, it is fair to assume that $1 \leq \omega_K (\xi) < 3$, where $\xi = z_1 / z_2 \in \C'$ is the slope of the start{i}ng point $\mathbf{z}$. For $\omega > \omega_K (\xi) \geq 1$, the first term in the sum on the right side of the inequality \eqref{E:errorirrfinal} dominates over the second, for all $k$'s large enough. Also, for $\omega > \omega_K (\xi),\ | q_k | \leq | q_{k - 1} |^{\omega}$ for $k \gg 0$ where we remind the reader that $p_k / q_k$'s are convergents to $\xi$ coming from a cont{i}nued fract{i}on expansion algorithm. Furthermore, under the condit{i}on~\eqref{E:jandk} and for $\omega < 3$, the second term in the upper bound for $| \gamma |$ is much smaller than the first implying that we have the existence of a $\gamma \in \Gamma$ which sat{i}sf{i}es
    \begin{equation}\label{E:gammaerror}
      | \gamma | \ll_{\mathbf{y}, \mathbf{z}, K} | q_k |^{\omega^{r_1 - 1} + 1}, \textrm{ and } | \gamma \mathbf{z} - \mathbf{y} | \ll_{\mathbf{y}, \mathbf{z}, K} \frac{1}{| q_k |^{\frac{1}{3 \omega} + \frac{4}{3} - \omega^{r_1 - 1}}}.
    \end{equation}
    The preceding lemma also tells us that $| \gamma | \gg | q_k q_{k - 1} |$ for the choice of $j$ and $k$ according to Eq.~\eqref{E:jandk}. This ensures the existence of inf{i}nitely many matrices $\gamma \in \SLtwoOK$ sat{i}sfying the above system of inequalit{i{es.\par
    As hinted before, we could have always started with an $\omega$ close enough to $1$ so that the exponent $\frac{1}{3 \omega} + \frac{4}{3} - \omega^{r_1 - 1}$ which expresses itself in the upper bound for $\gamma \mathbf{z} - \mathbf{y}$ in \eqref{E:gammaerror} is posit{i}ve. For such an $\omega > \omega_K (\xi)$, we therefore have that
    \begin{equation}
      \mu_{\Gamma} ( \mathbf{z}, \mathbf{y} ) \geq \frac{1 + 4 \omega - 3 \omega^{r_1}}{3 \omega ( \omega^{r_1 - 1} + 1 )}.
    \end{equation}
    In the limit $\omega_K (\xi) \leftarrow \omega$ from the right, and the generic value of the former being $1$, we get
    \begin{proposition}\label{P:irr}
      For all $\mathbf{z} \in \C^2 \setminus \{ \mathbf{0} \}$ having slope $\xi$ with $K$\,-\,irrat{i}onality measure $\omega_K (\xi) = 1$ and for all $\mathbf{y} \in \C^2$ with slope $y \in \C \setminus K$,
      \[
        \mu_{\Gamma} ( \mathbf{z}, \mathbf{y} ) \geq 1/3.
      \]
    \end{proposition}
    \noindent We now calculate lower bounds for $\hat{\mu}_{\Gamma}$. For almost all target points $\mathbf{y} \in \C^2 \setminus \{ \mathbf{0} \}$ with the slope $y$ belonging to $\C'$, we show a mildly stronger result than Prop.~\ref{P:irr} to be true, i.~e.,
    \begin{equation}
      \hat{\mu}_{\Gamma} ( \mathbf{z}, \mathbf{y} ) \geq 1/3
    \end{equation}
    for almost all pairs $( \mathbf{z}, \mathbf{y} ) \in \C^2 \times \C^2$. In proving this, the $K$\,-\,irrat{i}onality measure $\omega_K (y)$ associated with $\mathbf{y}$ is used as an auxiliary tool. The result below is the same as \cite[Lemma~6]{LN12} mutat{i}s mutandis and the proof is omit{t}ed.
    \begin{lemma}\label{L:tau}
      Let $\omega > \omega_K (\xi)\ ( \geq 1 )$ and def{i}ne
      \begin{equation}
        \tau := \frac{\omega_K (y)}{2 \omega_K (y) + 1} \omega^{r_1 - 1}.
      \end{equation}
      Given any $\eps > 0$ and $k_0 = k_0 (\eps) \in \N$, there exists $\gamma \in \Gamma$ such that
      \begin{equation}
        | \gamma | \ll_{\mathbf{y}, \mathbf{z}, K} | q_k |^{1 + \omega^{r_1 - 1}}\quad \textrm{and}\quad | \gamma \mathbf{z} - \mathbf{y} | \leq | q_k |^{\tau - 1 + \eps}
      \end{equation}
      for all $k > k_0$.
    \end{lemma}
    In addit{i}on to our assumpt{i}on that $\frac{1}{3 \omega} + \frac{4}{3} - \omega^{r_1 - 1}$ is posit{i}ve, we now also suppose an extra condit{i}on that $\omega < 2^{1 / ( r_1 - 1 )}$. This is to ensure that the quant{i}ty $\tau$ def{i}ned in Lemma~\ref{L:tau} remains less than $1$. Next, we restrict to $\eps$ small enough so that for given $\tau$, the exponent $1 - \tau - \eps$ whom we shall meet soon is greater than $0$. After this, our recourse is the old but very helpful idea of sandwiching (also used in \citeauthor{LN12}~\cite{LN12}) which given any suf{f}iciently large real positive number $T$, picks a $k$ large enough in terms of $T$ so that
    \begin{equation}
      C | q_k |^{1 + \omega^{r_1 - 1}} \leq T < C | q_{k + 1} |^{1 + \omega^{r_1 - 1}},
    \end{equation}
    where $C$ is the hidden constant in the upper bound for $| \gamma |$ given by Lemma~\ref{L:tau}. Such a choice of $k$ will mean that both
    \begin{equation}
      | \gamma | \leq T,\quad \textrm{and}\quad | \gamma \mathbf{z} - \mathbf{y} | \leq \frac{1}{| q_k |^{1 - \tau - \eps}} \leq \frac{1}{T^{(1 - \tau - \eps)/(\omega + \omega^{r_1})}}
    \end{equation}
    hold simultaneously, giving us a lower bound
    \begin{equation}
      \hat{\mu}_{\Gamma} ( \mathbf{z}, \mathbf{y} ) \geq \frac{1 - \tau - \eps}{\omega + \omega^{r_1}}.
    \end{equation}
    As the bound obtained is true for all suf{f}iciently small $\eps > 0$ and $\omega > \omega_K (\xi)$, in the limit $\eps \rightarrow 0_+$ and $\omega \rightarrow \omega_K (\xi)_+$,
    \begin{equation}
      \hat{\mu}_{\Gamma} ( \mathbf{z}, \mathbf{y} ) \geq \frac{\left( 2 - ( \omega_K (\xi) )^{r_1 - 1} \right) \omega_K (y) + 1 }{\left( 2 \omega_K (y) + 1 \right) \left( \omega_K (\xi) + ( \omega_K (\xi) )^{r_1} \right)}
    \end{equation}
    when the starting point $\mathbf{z} \in \C^2$ has slope whose $K$\,-\,irrat{i}onality measure $\omega_K (\xi)$ is very close to that of any generic point in the complex plane. Since we are only concerned with generic pairs $( \mathbf{z}, \mathbf{y} ) \in \C^2 \times \C^2$, we may as well take both $\omega_K (\xi)$ and $\omega_K (y)$ to be equal to $1$ whereby $\hat{\mu}_{\Gamma}$ comes out to be at least $1/3$.\par
    
    From the literature, we ment{i}on results of \citeauthor{Pol11}~\cite{Pol11} which is about calculat{i}ng the error term in the equidistribu{t}ion sum associated with the linear act{i}on of cocompact lat{t}ices $\Gamma \subset \SLtwoC$ on $\C^2$. The bounds for the generic value of Diophant{i}ne exponents then fall out as a corollary. A work in the same spirit for \SLtwoZ\ act{i}on on $\R^2 \setminus \{ \mathbf{0} \}$ was carried out by \citeauthor{MW12}~\cite{MW12} with much weaker est{i}mates than those of \citeauthor{LN12}~\cite{LN12}. Applicable in a broader framework, the machinery of \citeauthor*{GGN3}~\cite{GGN3, GGN4} is vastly superior and gives the values of exponents for an array of lat{t}{i}ce act{i}ons on homogeneous variet{i}es of connected almost simple, semisimple algebraic groups (see, in particular~\cite{GGN3, GGN4}). However, for them too, the lower bound for the uniform exponent $\hat{\mu}_{\Gamma}$ as in Def.~\ref{D:dexpunf} for any given $\mathbf{z}$ with dense $\Gamma$\,-\,orbits in the complex plane and a generic target point $\mathbf{y} \in \C^2$ is of{f} by a factor of $2$ compared to ours, as has been told by the authors in a personal communicat{i}on. Nevertheless, it is one of the papers in this series that we look at next in our search for upper bounds on $\mu_{\Gamma}$ and $\hat{\mu}_{\Gamma}$.\par
    
    Let us turn our at{t}ent{i}on to Theorem~3.1 of~\cite{GGN3}. The punctured plane $\C^2 \setminus \{ \mathbf{0} \}$ can be realized as the special linear group $\SLtwoC$ quot{i}ented by the closed upper unipotent subgroup $H$. The  non-uniform lat{t}ice $\Gamma = \SLtwoOK$ acts ergodically on $G / H$ and we verify that the hypothesis of the theorem is valid in this scenario. In the terminology of \citeauthor{GGN3}, the coarse volume growth exponent $a$ for the upper unipotent group $H \subset \SLtwoC = G$ and the lower local dimension $d'$ of the homogeneous space $G / H \approx \C^2 \setminus \{ \mathbf{0} \}$ equal $2$ and $4$, respect{i}vely. We then have that for any $\mathbf{z} \in \C^2$ with a dense $\Gamma$-orbit and almost all $\mathbf{y} \in \C^2$, the inverse
    \begin{equation}
      \kappa ( \mathbf{z}, \mathbf{y} ) := \frac{1}{\hat{\mu}_{\Gamma} ( \mathbf{z}, \mathbf{y} )} \geq 2
    \end{equation}
    which is the same as saying that $\hat{\mu}_{\Gamma} ( \mathbf{z}, \mathbf{y} ) \leq 1/2$ for all $\mathbf{z} \in \C^2$ with slope $\xi \in \C'$ and $\mathbf{y}$ belonging to a full Lebesgue measure subset (depending on $\mathbf{z}$) of the complex plane. The same proof can be modif{i}ed to replace $\hat{\mu}_{\Gamma}$ with $\mu_{\Gamma}$. To see this, one should apply Borel-Cantelli as soon as we have the estimates for the number of lat{t}{i}ce elements $\gamma \in \Gamma \cap G_t$  of bounded size $e^t$ and such that $\gamma \mathbf{z}$ lies within a unit distance of the target point $\mathbf{y}$. This is given by $e^{2 t + \eps}$ upto a constant multiple depending on $\eps$ alone. In summary,
    \begin{proposition}\label{P:generic}
      For all pairs $( \mathbf{z}, \mathbf{y} ) \in \C^2 \times \C^2$ with the slopes $\xi$ of $\mathbf{z}$ and $y$ of $\mathbf{y}$ both having $K$\,-\,irrat{i}onality measure equal to $1$, we have
      \begin{equation}
        1/3 \leq \hat{\mu}_{\Gamma} ( \mathbf{z}, \mathbf{y} ) \leq \mu_{\Gamma} ( \mathbf{z}, \mathbf{y} ), \textrm{ and}
      \end{equation}
      for all $\mathbf{z} \in \C^2$ with slope $\xi \in \C'$ and almost all (depending on $\mathbf{z}$) target point $\mathbf{y}$, the upper bound
      \begin{equation}
        \hat{\mu}_{\Gamma} ( \mathbf{z}, \mathbf{y} ) \leq \mu_{\Gamma} ( \mathbf{z}, \mathbf{y} ) \leq 1/2.
      \end{equation}
    \end{proposition}

  \subsection{\protect{Target point with $K$\,-\,rat{i}onal slope}}\label{SS:rat}
    The task of comput{i}ng exponents is much easier when $\mathbf{y}$ has slope $y = y_1 / y_2 = a / b \in K$, where $a, b \in \O_K$ and $| \gcd_{\O_K} ( a, b ) | = 1$. Without any loss of generality, assume as before that $\max \{ 1, |a| \} \leq |b|$. The column vector $( a, b )^t$ is taken to be the first column of our matrix $N$ as the fraction $a/b$ is the best approximat{i}on to $y$ by $K$\,-\,rat{i}onal points. Af{t}er this step, the second column can be chosen to be some $(a', b')^t \in \O_K^2$ such that $ab' - a'b = 1$ and $|b'| \leq |b|$. This is possible because of the fact that $\O_{\Q (\sqrt{-d})}$ is a Euclidean domain for $d = 1, 2, 3, 7$ and $11$~\cite{Mot49} and then, it is clear that $|N| \asymp |b|$. As in Sec.~\ref{SS:generic}, let $\omega > \omega_K (\xi)\ (\geq 1)$ where $\xi = z_1 / z_2 \in \C'$. If necessary, we will take $\omega$ to be very close to $\omega_K (\xi)$.
    \begin{lemma}[\protect{cf.~\cite[Lemma~5]{LN12}}]\label{L:gammarat}
      Let $k \in \N$ be large enough. Given $\mathbf{y} \in \C^2$ with slope $y \in K$, there exists some $\ell \in \O_K$ and $\gamma = N U^{\ell} M_k \in \Gamma$ sat{i}sfying
      \begin{align}
        \abs{ q_k q_{k - 1} } \ll_{\mathbf{y}, \mathbf{z}, K} \abs{ \gamma } \ll_{\mathbf{y}, \mathbf{z}, K} &\abs{q_{k - 1}} \abs{ q_k }^{\omega^{r_1 - 1}} \textrm{ and also,}\notag\\
        \abs{ \gamma \mathbf{z} - \mathbf{y} } \ll_{K} &\abs{\frac{b z_2}{q_k}}.
      \end{align}
    \end{lemma}
    \begin{proof}
      Here, the quant{i}t{i}es $\delta$ and $\delta'$ def{i}ned in Lemma~\ref{L:gammay} equal $by - a = 0$ and $b'y - a' = 1/b$, respect{i}vely and thereby,
      \begin{equation}
        | \Lambda_1 - y \Lambda_2 | \ll_{K} \left| \frac{z_2}{b q_k} \right|.
      \end{equation}
      For the same reason, $b$ replaces $s_j$ in the inequality \eqref{E:lambda2bound} and, thereaf{t}er, the triangle inequality gives the required bound on $\abs{ \gamma \mathbf{z} - \mathbf{y} }$. The same change made to Eq.~\eqref{E:ellrange} gives us
      \begin{equation}
        \abs{ \gamma } \ll \abs{ \ell b q_{k - 1} } + \abs{ b q_k } \ll_{\mathbf{y}, \mathbf{z}, K} \abs{q_{k - 1}} \abs{ q_k }^{\omega^{r_1 - 1}} + \abs{q_{k - 1}} + \abs{q_k} \ll \abs{q_{k - 1}} \abs{ q_k }^{\omega^{r_1 - 1}}
      \end{equation}
      in conjunct{i}on with Lemma~\ref{L:gammabound}. On the other hand,
      \begin{equation}
        \abs{ \gamma } \gg_{\mathbf{z}} \abs{ \ell q_{k - 1} } - \abs{q_k} \gg_{\mathbf{y}, \mathbf{z}, K} \abs{ q_k q_{k - 1} }
      \end{equation}
      while again looking at the modif{i}ed Eq.~\eqref{E:ellrange}, as $\abs{ q_{k - 1} } \gg \abs{ b }$ for all large enough $k$'s.
    \end{proof}
    \noindent From the bounds on the size of $\gamma$ in the above lemma, we get
    \begin{equation}
      \abs{ q_k } \ll_{\mathbf{y}, \mathbf{z}, K} \abs{ \gamma } \ll_{\mathbf{y}, \mathbf{z}, K} \abs{ q_k }^{\omega^{r_1 - 1} + 1}.
    \end{equation}
    The second part in this inequality gives us
    \begin{equation}
      \abs{ \frac{1}{q_k} } \ll_{\mathbf{y}, \mathbf{z}, K}
      \abs{ \gamma }^{- \frac{1}{\omega^{r_1 - 1} + 1}}
    \end{equation}
    which in turn implies that
    \begin{equation}
      \abs{ \gamma \mathbf{z} - \mathbf{y} } \ll_{\mathbf{y}, \mathbf{z}, K}
      \abs{ \gamma }^{- \frac{1}{\omega^{r_1 - 1} + 1}}
    \end{equation}
    for inf{i}nitely many $\gamma \in \SLtwoOK$. This means that the Diophant{i}ne exponent $\mu_{\Gamma} ( \mathbf{z}, \mathbf{y} ) \geq ( \omega^{r_1 - 1} + 1 )^{-1}$ for all     $\mathbf{y}$ with slope $y \in K$ and $\omega > \omega_K (\xi)$. Taking the limit $\omega_K (\xi) \leftarrow \omega$ from the right, we conclude that for any start{i}ng point $\mathbf{z} \in \C^2 \setminus \{ \mathbf{0} \}$ whose slope $\xi$ has $K$\,-\,irrat{i}onality measure $\omega_K (\xi)$ and any target point $\mathbf{y}$ with ``$K$\,-\,rat{i}onal slope'',
    \begin{equation}\label{E:murat}
      \mu_{\Gamma} ( \mathbf{z}, \mathbf{y} ) \geq \frac{1}{\left( \omega_K (\xi) \right)^{r_1 - 1} + 1}
    \end{equation}
    Next on our agenda is a lower bound for $\hat{\mu}_{\Gamma}$ when $\mathbf{y}$ has a $K$\,-\,rat{i}onal slope. This value is in general lower than that for $\mu_{\Gamma}$ above, but equals the same for almost every $\mathbf{z}$. Given $T \gg 0$ and $\omega > \omega_K (\xi)$, we choose $k$ as
    \begin{align}
      \abs{q_{k - 1}} \abs{q_k}^{\omega^{r_1 - 1}} \leq T &< \abs{q_k} \abs{q_{k + 1}}^{\omega^{r_1 - 1}} \textrm{ whereby}\notag\\
      \abs{\gamma} \ll_{\mathbf{y}, \mathbf{z}, K} \abs{q_{k - 1}} \abs{q_k}^{\omega^{r_1 - 1}} &\leq T, \textrm{ and } T \leq \abs{q_k}^{1 + \omega^{r_1}}
    \end{align}
    for $\gamma$ given to us by Lemma~\ref{L:gammarat}. No more input, apart from repeating the same set of arguments, is required to now deduce that
    \begin{equation}\label{E:hatmurat}
      \hat{\mu}_{\Gamma} ( \mathbf{z}, \mathbf{y} ) \geq \frac{1}{\left( \omega_K (\xi) \right)^{r_1} + 1}.
    \end{equation}
    Borrowing an idea from \citeauthor{LN12b}~\cite{LN12b}, we get the following transference result
    \begin{proposition}\label{P:minkowski}
      Let $\xi, y \in \C$. Then, the exponents for inhomogeneous approximat{i}on by $\O_K$\,-\,integers
      \[
        \hat{\omega}_K ( \xi, y ) \geq \frac{1}{\left( \omega_K (\xi) \right)^{r_1} + 1}\quad \textrm{and}\quad \omega_K ( \xi, y ) \geq \frac{1}{\left( \omega_K (\xi) \right)^{r_1 - 1} + 1}.
      \]
    \end{proposition}
    \begin{proof}
      In the above observations, let $\mathbf{z} = ( \xi, 1 )^t$ and $\mathbf{y} = ( y, y )^t$. Either of the two rows of the various matrix solut{i}ons $\{ \gamma_i \} \subset \SLtwoOK$ thus obtained will do the job.
    \end{proof}
    As a special case when $\omega_K (\xi) = 1$, we obtain that for all $\eps > 0$, there exists a $T_0 > 0$ and for all $T > T_0$, we have a pair $( q, p ) \in \O_K^2$ for which
    \begin{equation}\label{E:minkowski}
       \abs{ q \xi + p - y } < \frac{1}{T^{1/2 - \eps}}\quad \textrm{and}\quad \max \{ \abs{p}, \abs{q} \} \leq T.
    \end{equation}
    The lemma writ{t}en below helps us to obtain an upper bound for the Diophant{i}ne exponent $\mu_{\Gamma} ( \mathbf{z}, \mathbf{y} )$ when the starting point has dense \SLtwoOK-orbit in $\C^2$ and the target point has a $K$\,-\,rat{i}onal slope. The method used is \citeauthor{LN12}'s factorizat{i}on technique \cite[Theorem~4]{LN12} to break down any candidate matrix $\gamma \in \Gamma$ in terms of well-known ent{i}t{i}es like $N$ and $M_k$ in order to be able to say something about the size of $\gamma \mathbf{z} - \mathbf{y}$ and of the various components appearing in between. As $\hat{\mu}_{\Gamma} \leq \mu_{\Gamma}$, this will trivially give us an upper bound for $\hat{\mu}_{\Gamma} ( \mathbf{z}, \mathbf{y} )$.
    \begin{lemma}
      Let $\mathbf{z} \in \C^2 \setminus \{ \mathbf{0} \}$ have a slope $\xi \in \C'$ and $\mathbf{y}$ be a fixed target point with slope $y = a/b \in K$ as introduced in the beginning of the sect{i}on. For all $k$ large enough and $\gamma \in \Gamma$ such that
      \[
        \abs{\gamma} \leq \frac{1}{3 C_1}\abs{\frac{y_2}{z_2}} \abs{ q_k q_{k + 1}},\quad \textrm{we must have}\quad \abs{\gamma \mathbf{z} - \mathbf{y} } \geq \abs{\frac{z_2}{3b}} \frac{1}{\abs{q_k}}
      \]
    \end{lemma}
    \noindent Here, $C_1$ refers to the constant discussed in Eq.~\eqref{E:uerror}, dist{i}lled from \citeauthor{Dan15}'s cont{i}nued fract{i}on theory for complex numbers in terms of $\O_K$\,-\,integers.
    \begin{proof}
      Assume, if possible, that for some $\gamma = \begin{pmatrix} v_1 & u_1\\v_2 & u_2 \end{pmatrix}$ as above, the vector $\gamma \mathbf{z} - \mathbf{y} = \begin{pmatrix} \Lambda_1\\ \Lambda_2 \end{pmatrix}$ has supremum norm strictly less than $\abs{z_2} \abs{3b q_k}^{-1}$. Without loss of generality, we may suppose that $|a| \leq |b|$ because of the matrix $J$ from Eq.~\eqref{E:J}. Given the complex number $a/b$ with $a, b \in \O_K$ and $| \gcd_{\O_K} ( a, b ) | = 1$, we take $N = \begin{pmatrix} a & a'\\ b & b' \end{pmatrix} \in \Gamma$ with $|b| \leq |N| < 2 |b|$. As in \cite{LN12}, let
      \begin{equation}
        \gamma' := N^{-1} \gamma = \begin{pmatrix} v'_1 & u'_1\\v'_2 & u'_2 \end{pmatrix}.
      \end{equation}
      Since $b' y_1 - a' y_2 = y_2 / b$, here too we get that
      \begin{equation}\label{E:gammaprime}
        \gamma'
        = \begin{pmatrix}
            \frac{b' ( v_1 y_2 - v_2 y_1)}{y_2} + \frac{v_2}{b} & \frac{b' ( u_1 y_2 - u_2 y_1)}{y_2} + \frac{u_2}{b}\\
            - \frac{b ( v_1 y_2 - v_2 y_1)}{y_2}                & - \frac{b ( u_1 y_2 - u_2 y_1)}{y_2}
          \end{pmatrix}
      \end{equation}
      af{t}er adding and subtract{i}ng equal quant{i}ties to both the entries in the first row. Also,
      \begin{equation}\label{E:gammaprimez}
        \begin{pmatrix}
          z_2 ( v'_1 \xi + u'_1 )\\
          z_2 ( v'_2 \xi + u'_2 )
        \end{pmatrix} =
        \gamma' \mathbf{z} 
        = N^{-1}
        \left(
        \mathbf{y} +
        \begin{pmatrix}
          \Lambda_1\\
          \Lambda_2
        \end{pmatrix} \right) =
        \begin{pmatrix}
          \frac{y_2}{b} + b' \Lambda_1 - a' \Lambda_2\\
          -b \Lambda_2 + a \Lambda_2
        \end{pmatrix}.
      \end{equation}
      Next, the determinant
      \begin{align}
        v_1 y_2 - v_2 y_1 &=
        \begin{vmatrix}
          v_1 & y_1\\
          v_2 & y_2 
        \end{vmatrix} = 
        \begin{vmatrix}
          \begin{matrix}
            v_1\\
            v_2
          \end{matrix} & \gamma \mathbf{z} -
          \begin{pmatrix}
            \Lambda_1\\
            \Lambda_2
          \end{pmatrix}
        \end{vmatrix} = z_2 -
        \begin{vmatrix}
          v_1 & \Lambda_1\\
          v_2 & \Lambda_2
        \end{vmatrix}, \textrm{ whereby}\notag\\
        \abs{ v_1 y_2 - v_2 y_1 } &\leq \abs{z_2} + 2 \abs{\gamma} \max \{ \abs{\Lambda_1}, \abs{\Lambda_2} \}\label{E:det}\\
                                  &\leq \abs{z_2} + \abs{ \frac{2 y_2}{9 C_1 b} q_{k + 1}} \leq \abs{ \frac{y_2}{4 C_1 b} } \abs{q_{k + 1}}\notag
      \end{align}
      for all $k$ such that $|q_k| > 36 C_1 \abs{ b z_2 } / \abs{ y_2 }$. Combining the last three Eqs.~\eqref{E:gammaprime}, \eqref{E:gammaprimez} and \eqref{E:det}, we conclude that
      \begin{align}
        \abs{v'_2}              &= \abs{\frac{b}{y_2} ( v_1 y_2 - v_2 y_1 )} < \frac{\abs{ q_{k + 1} }}{4 C_1}, \textrm{ and}\notag\\
        \abs{ v'_2 \xi + u'_2 } &= \frac{1}{\abs{z_2}} \abs{ -b \Lambda_1 + a \Lambda_2 } \leq 2 \abs{\frac{b}{z_2}} \max \{ \abs{\Lambda_1}, \abs{\Lambda_2} \} < \frac{2}{3} \cdot \abs{q_k}^{-1}
      \end{align}
      as $\abs{a} \leq \abs{b}$. Now, consider the \SLtwoOK\ matrix $g$ def{i}ned to be $N^{-1} \gamma M_k^{-1}$. Then, $g = \begin{pmatrix} * & *\\ * & v'_2 p_k + q_k u'_2 \end{pmatrix}$, and the lower left entry has size
      \begin{align}
        \abs{v'_2 p_k + q_k u'_2} &= \abs{- v'_2 (q_k \xi - p_k ) + q_k ( v'_2 \xi + u'_2 ) }\notag\\
                                  &\leq \abs{\frac{C_1 v'_2}{q_{k + 1}}} + \abs{q_k} \abs{v'_2 \xi + u'_2} < \frac{1}{4} + \frac{2}{3} < 1.
      \end{align}
      As the ring of integers $\O_K$ was taken to be discrete, it has no non-zero element whose Euclidean norm is less than one. This means $g$ has to be equal to $\begin{pmatrix} m & \zeta\\ - \zeta^{-1} & 0 \end{pmatrix}$ for some $m \in \Z$ and $\zeta \in \O_K^*$ with $\abs{\zeta} = 1$. Thereaf{t}er, Eq.~\eqref{E:gammaprimez} tells us that the vector
      \begin{equation}
        \gamma' \mathbf{z} =
        \begin{pmatrix}
          \frac{y_2}{b} + b' \Lambda_1 - a' \Lambda_2\\
          -b \Lambda_2 + a \Lambda_2
        \end{pmatrix} = g M_k \mathbf{z} = z_2
        \begin{pmatrix}
          m \epsilon_k + (-1)^{k - 1} \zeta \epsilon_{k - 1}\\
          - \zeta^{-1} \epsilon_k
        \end{pmatrix}.
      \end{equation}
      We concentrate on the entry in the f{i}rst coordinate to get
      \begin{align}
        \abs{\frac{y_2}{b}} - \frac{4}{3} \abs{\frac{z_2}{q_k}} \leq \abs{\frac{y_2}{b} + b' \Lambda_1 - a' \Lambda_2} &= \abs{z_2 }\abs{m \epsilon_k + (-1)^{k - 1} \zeta \epsilon_{k - 1}}\notag\\
        &\leq C_1 \abs{z_2} \left( \abs{\frac{m}{q_{k + 1}}} + \abs{\frac{1}{q_k}} \right)
      \end{align}
      which gives a lower bound on $\abs{m}$,
      \begin{equation}
        \abs{m} \geq \abs{\frac{101 y_2}{108 C_1 b z_2}} \abs{q_{k + 1}} > 33
      \end{equation}
      as $C_1 > 1$ and we recall that $|q_{k + 1}| > |q_k| > 36 C_1 \abs{ b z_2 } / \abs{ y_2 }$. We now have a decomposit{i}on for the matrix $\gamma$ as $\gamma = N g M_k$ which helps us to get a handle on its size. To be precise,
      \begin{equation}
        \gamma = \pm
        \begin{pmatrix}
          a & a'\\
          b & b'
        \end{pmatrix}
        \begin{pmatrix}
          m & -1\\
          1 & 0
        \end{pmatrix}
        \begin{pmatrix}
          q_k                    & -p_k\\
          (-1)^{k - 1} q_{k - 1} & (-1)^k p_{k - 1}
        \end{pmatrix}
      \end{equation}
      so that $\abs{\gamma} \geq \abs{b} ( \abs{m} - 2 ) \abs{q_k}$ by the triangle inequality. As we have argued that $\abs{m}$ should be greater than $33$ for all suitably large $k$'s, the quant{i}ty $\abs{m} - 2$ should be strictly bigger than $31 \abs{m} / 33$. We, therefore, deduce that
      \begin{equation}
        \abs{\gamma} > \frac{31}{33} \abs{b m q_k} \geq \frac{4}{5C_1} \abs{\frac{y_2}{z_2}} \abs{ q_k q_{k + 1} },
      \end{equation}
      but the hypothesis was $\abs{\gamma} \leq \frac{\abs{y_2}}{3 C_1 \abs{z_2}} \abs{ q_k q_{k + 1}}$, a contradict{i}on.
    \end{proof}
    For any $\gamma \in \Gamma$ with $\abs{\gamma}$ suf{f}iciently large, pick $k$ as
    \begin{equation}
      \frac{1}{3 C_1}\abs{\frac{y_2}{z_2}} \abs{ q_{k - 1} q_k} < \abs{\gamma} \leq \frac{1}{3 C_1}\abs{\frac{y_2}{z_2}} \abs{ q_k q_{k + 1}},
    \end{equation}
    and for the case when $\omega_K (\xi)$ is f{i}nite, choose any real number $\omega > \omega_K (\xi)$ so that eventually we have $| q_{k - 1} | \geq \abs{q_k}^{1 / \omega}$. Consequently,
    \begin{equation}
      \abs{ \gamma \mathbf{z} - \mathbf{y} } \geq \abs{\frac{z_2}{3b}} \frac{1}{\abs{q_k}} \gg_{ \mathbf{y}, \mathbf{z} } \frac{1}{\abs{\gamma}^{\omega / ( \omega + 1 ) }}
    \end{equation}
    and let{t}ing $\omega$ approach $\omega_K (\xi)$ from the right, we get that
    \begin{equation}\label{E:muratub}
      \mu ( \mathbf{z}, \mathbf{y} ) \leq \frac{\omega_K (\xi)}{\omega_K (\xi) + 1}.
    \end{equation}
    The statement is also true for $\omega_K (\xi) = \infty$, as can be easily checked. However, as discussed in Sect{i}on~\ref{S:intro}, the $K$\,-\,irrat{i}onality measure $\omega_K (\xi)$ equals $1$ for Lebesgue\,-\,almost all $\xi \in \C$. Therefore, for a full measure subset of $\C^2 \setminus \{ \mathbf{0} \}$, Eqs.~\eqref{E:hatmurat} and \eqref{E:muratub} together give
    \begin{proposition}\label{P:rat}
      For all $\mathbf{z} \in \C^2 \setminus \{ \mathbf{0} \}$ such that the slope $\xi$ has $K$\,-\,irrat{i}onality measure to be $1$, and for all $\mathbf{y} \in \C^2$ with slope $y \in K$,
      \[
        \mu_{\Gamma} ( \mathbf{z}, \mathbf{y} ) = \hat{\mu}_{\Gamma} ( \mathbf{z}, \mathbf{y} ) = 1/2.
      \]
    \end{proposition}
    This completes the proof of Theorem~\ref{Th:main}.

\bibliographystyle{plainnat}
\bibliography{gaussian_lattice}
\end{document}